# Efficiency Evaluation of Banks with Many Branches using a Heuristic Framework and Dynamic Data Envelopment Optimization Approach: A Real Case Study

**Vahid Kayvanfar[1] , Hamed Baziyad[2], Shaya Sheikh[3], Frank Werner[4]**

[1]*Division of Engineering Management and Decision Sciences, College of Science and Engineering, Hamad Bin Khalifa University, Doha, Qatar*
*email: valikayvanfar@hbku.edu.qa*

[2]*Department of Information Technology Engineering, Faculty of Industrial and Systems Engineering, Tarbiat Modares University, Tehran, Iran*
*email: hamed_baziyad@modares.ac.ir*

[3] *School of Management, New York Institute of Technology, 1855 Broadway, New York, NY 10023, United States*
*email: ssheik11@nyit.edu*

[4] *Faculty of Mathematics, Otto-von-Guericke University Magdeburg, P.O. Box 4120, 39016 Magdeburg, Germany*
*email: frank.werner@ovgu.de*

## Abstract

Evaluating the efficiency of organizations and branches within an organization is a challenging issue for managers. Evaluation criteria allow organizations to rank their internal units, identify their position concerning their competitors, and implement strategies for improvement and development purposes. Among the methods that have been applied in the evaluation of bank branches, non-parametric methods have captured the attention of researchers in recent years. One of the most widely used non-parametric methods is the data envelopment analysis (DEA) which leads to promising results. However, the static DEA approaches do not consider the time in the model. Therefore, this paper uses a dynamic DEA (DDEA) method to evaluate the branches of a private Iranian bank over three years (2017-2019). The results are then compared with static DEA. After ranking the branches, they are clustered using the K-means method. Finally, a comprehensive sensitivity analysis approach is introduced to help the managers to decide about changing variables to shift a branch from one cluster to a more efficient one.

**Keywords**: Dynamic data envelopment analysis (DDEA); Branch performance evaluation; Efficiency; Clustering; Banking industry; Big data.

## 1. Introduction

The structural design for the identification, evaluation, and ultimately improving operations in branches deems necessary. Financial markets play an important role in the economy of countries, and among them, banks have a special role as the most important element of the financial market. In this regard, recognizing the types of efficiency and measuring the impact of factors and policies on their efficiency has drawn the attention of researchers. Such studies are of critical importance both from a microeconomics as well as a macroeconomics perspective (Ariff & Can, 2008; George Assaf et al., 2011; Chan et al., 2015). From the microeconomic lens, efficiency intensifies competition and improvements in institutional and regulatory frameworks. From the macroeconomic perspective, the banking system's efficiency impacts the cost of financial intermediation and the stability of financial markets. Any improvement in bank performance signals a better allocation of financial resources and, as a result, an increase in the investment. For the past two decades, there has been a steady increase in the breadth and depth of studies focused on the banking system's efficiency worldwide (Tsolas et al., 2020; Preeti & Roy, 2020; Fukuyama et al., 2020).

A performance evaluation is carried out in either parametric (econometrics) or non-parametric (linear/mathematical programming) format, each of which has several strengths and weaknesses in the literature. At the same time, there is no explicit consensus on the superiority of one over the other approach. In recent years, researchers have recommended both methods to analyze the same data set, demonstrate the robustness and awareness of the results, and provide more reliable information and meaningful implications for policymakers and managers. According to Bauer et al., (1998), instead of forming a consensus on the superiority of one method, it is enough to introduce a set of compatibility conditions that the results of the proposed model must meet to make it worthwhile for regulators.

Researchers in several countries have proposed a wide range of methodologies for evaluating the performance of banks; see Shahzad Virk et al., (2022), Chen et al., (2022), and Ben Lahouel et al., (2022), (Gkougkousi et al., 2022). Financial and monetary enterprises, especially commercial banks, have a direct and significant impact on the economic development of the countries. As a result, the existence of a model that provides feedback and improves different branches' performance is critical. Numerous institutions around the world rank companies. Depending on the institution's mission, different quantitative and qualitative indicators are used for these rankings. *Standard & Poor*, *Moody's, Fortune*, *Business Week*, and *Financial Times* are some of these institutions. In Iran, the *Industrial Management Organization* and the *Stock Exchange and Securities Organization* rank the top companies. The *Industrial Management Organization* annually presents a list of the top 100 companies. The Exchange and Securities Organization also ranks Iranian stock exchange companies with the aim of facilitating the selection of suitable companies for an investment.

The DEA methodology has been used to measure the relative efficiency of decision-making units (DMUs), and it is used as one of the most applied tools in decision-making, particularly for a bank branches' performance analysis (Balak et al., 2021; Chu et al., 2022; Omrani et al., 2022). This method measures and ranks the performance of existing companies based on inputs and outputs. The DMUs are decision-making units with multiple inputs and outputs (Hu et al., 2009). DEA is based on a series of linear programming optimizations, also called the non-parametric method. In this method, an efficient boundary frontier is created from a series of points that are determined by linear programming. Either the constant or variable returns to scale assumption can be used to determine these points. Return-to-scale represents the relationship between changes in the inputs and outputs of a system. The linear programming method determines whether the DMU is on or outside the performance boundary. In this way, efficient and inefficient DMUs are separated from each other. Charnes et al., (1978) provided a model that measures the performance with multiple inputs and outputs (a.k.a. CCR model). CCR was first used in assessing the academic achievement of United States school students in 1976. Banker et al., (1984) embedded new concepts into DEA and called it the BCC model. The CCR model has a constant return-to-scale that tries to achieve the efficiency of a DMU by selecting the optimal weights for the input and output variables. Unlike the CCR model, the BCC model works based on the variable return-to-scale assumption. Using variable returns to scale makes it possible to provide a very accurate analysis by calculating the technical efficiency in terms of scale and management efficiency values. The same principles of the CCR model are used to build input and output-oriented models in the BCC model. In the input-oriented model, the efficiency increases with decreasing the inputs; in the output model, the efficiency increases with increasing the outputs (Shao et al., 2018; Kazemi Matin et al., 2021).

An advantage of DEA compared to other methods is that DEA is not sensitive to the unit of measurement. In this method, inputs and outputs can have different units. Another advantage is that DEA compares the units with each other. Despite multi-criteria decision analysis (MCDM) approaches, DEA models do not need experts' opinions for weighting the criteria because the optimum values of the weights are calculated by themselves, which is an appropriate tool for dealing with uncertainty (Filatovas et al., 2022; Yi-Chia et al., 2022). However, DEA models suffer from a big data analysis, which makes the feasibility area smaller and may lead to the infeasibility of the problem under investigation (Badiezadeh et al., 2018; Zhou et al., 2021; Zhu, 2022). Nowadays, big data generation has become challenging in different areas, such as text mining (Hosseini et al., 2021; Pourhatami et al., 2021), healthcare (Shirazi et al., 2020), demand forecasting (Zohdi et al., 2022), and the Internet of Things (Baziyad et al., 2022). Therefore, utilizing fewer data (Norouzi et al., 2022) or developing methods and approaches with the capability of handling a large amount of data (Anagnostopoulos et al., 2016) was changed into an exciting area. This paper proposes a heuristic framework, enabling DEA models to be run in big data environments. The proposed framework can also enable DEA models to deal with infeasible problems by changing them from infeasible to feasible ones. This process is executed by eliminating some constraints from the model.

A comprehensive overview of the previous studies of DEA-based bank performance analysis is provided in Table 1. Different inputs and outputs are provided in different countries as needed. In the case of inputs, variables such as total cost, interest, capital, profit and income, shareholders 'salaries, employees' salaries, number of employees, total assets and fixed assets, as well as deposits, are among the most frequently used items. Among them, a large margin uses items such as the *number of employees*, *staff costs*, *deposits*, *fixed assets*, and *total assets*. Some variables such as *income and profit*, *loans*, *earnings per share*, *investment*, *securities and bonds*, *deposits*, and *shareholder returns* are used more frequently.

Among the many methods that have been registered in the evaluation of bank branches, the interest of researchers in non-parametric methods has increased in recent years, which constitutes a major part of the evaluation of organizations. One of the most widely used non-parametric methods in this field is the use of DEA techniques (J. S. Liu et al., 2016). DEA can be divided into static and dynamic categories (von Geymueller, 2009). In the static state, the desired time intervals are not realized, and the research is performed at one point, which may lead to a misleading data analysis. In contrast, dynamic DEA generates more realistic results since it considers the changes that occur over time. In problems where DEA traditionally uses a performance appraisal, the time factor is not considered, despite having time as one of the most important performance factors. Eliminating time intervals could create errors in the output, and as a result, a wrong conclusion by managers becomes more likely. For this purpose, dynamic DEA is used to consider time intervals in evaluating the branches of a private bank in Iran over three years. The results are then compared with the classical static DEA. After ranking the branches, a clustering of the branches is performed using the K-means method. Finally, to present a situational analysis of the branches and help managers to plan for the future, the amount of change needed in each key indicator (to switch the cluster membership of each branch) is obtained using a comprehensive sensitivity analysis. The research goals can be summarized as follows:

- Investigating the efficiency of the branches of a private bank throughout the country for three years from 2017 to 2019.
- Implementing dynamic DEA in evaluating the bank branches for the first time.
- Ranking of the private bank branches during similar time intervals.

Table 1: Summary of DEA researches conducted by country and type of model

| Researcher and year | Goal | Country | Model | Input | Output | No. DMUs |
|---|---|---|---|---|---|---|
| (Kamarudin et al., 2015) | TE, PTE and SE review of 43 banks from 2007 to 201 | Gulf Cooperation Council (GCC) countries | CCR–1 and BCC–1 | Large assets and deposits | Loans and income | 43 |
| (Kwon & Lee, 2015) | Combining two experimental data analysis techniques to evaluate and predict the performance improvement | USA | Two stages: CCR-2 and Back Propagation Neural Network (BPNN) | Step I: Number of employees, equity and costs.<br>Stage II: deposits, loans and investments | Step I: Deposits, Loans and Investments.<br>Step II: Profit | 13 |
| (Yilmaz & Güneş, 2015) | Comparing the technical and general productivity of ordinary deposit banks and Islamic banks from 2007 to 2013 | Turkey | CCR-1 and BCC-1 | Deposits and capital | Loans, investments and income | 32 |
| (Stoica et al., 2015) | Analyzing the Internet banking services' performance | Romania | PCA-based DEA | Deposits, total deposits, remittances and Total operating costs, Number of employees, The value of owned equipment and software programs | Net total revenues, Daily "reach" average rate | 24 |
| (Wanke et al., 2016) | A new fuzzy DEA model to evaluate the bank performance in Mozambique from 2003 to 2011 | Mozambique | Fuzzy | Total costs (excluding personnel) and personnel costs | Total deposits, pre-tax income and total credit operations | 13 |
| (Kamarudin et al., 2016) | Effects of government on the revenue productivity of Islamic and conventional banks in the period 2007-2011 | GCC countries | Two-stage - BCC with panel regression | Staff costs and deposits | Income and loans | 47 |
| (Stewart et al., 2016) | The efficiency of the Vietnamese banking system from 1999 to 2009 | Vietnam | Two stages: CCR, BCC and bootstrap | Number of employees, deposits from other banks and customer deposits | Loans from customers, other loans and securities | 101 |
| (Kamarudin et al., 2017) | Productivity of Islamic banks by confirming specific factors of banks, industry | Southeast Asian countries | Two-stage - BCC with Malmquist and panel regressions | Staff costs, deposits and fixed assets | Investments and income | 29 |
| (Fukuyama & Matousek, 2017) | Evaluating the banks' network revenue performance | Japan | | Number of workers, Premises and real estate | Performing loans, Securities investment, Nonperforming loans | 172 |
| (X. Liu et al., 2020) | Finding efficiency gap among different commercial banks | China | Two-stage DEA | Fixed assets, Labor Operating expenses | Non-performing loan ratio, Loans Business income | 28 |
| (Tsolas et al., 2020) | Classifying the sampled branches of banks | Greece | Two-stage hybrid DEA- ANN | personnel expenses, rents and depreciation, operational expense | Net interest income, non-interest income, | 160 |
| (Preeti & Roy, 2020) | Operational performance measurement and prediction | India | DEA+ANN | Operating expenses and deposits | Investments, performing loans and advances, non-interest income | 39 |
| (Preeti et al., 2020) | Comparing the efficiency of private and public banks | India | Frontier-based radial DEA | Operating expenses, Assets and Deposits | Investment and Advances | 40 |
| (Fukuyama et al., 2020) | Analyzing the cost inefficiency levels of Turkish banks | Turkey | two-stage DEA | Total fixed assets, Number of employees | Various loans, Security investments | 26 |
| (Li, 2020) | Analyzing the performance of city commercial banks in China | China | Bootstrapped DEA | Fixed assets, Staff expenses, Interest expenses | Non-performing loans, Deposits, Performing loans, Non-interest income | 101 |
| (Lartey et al., 2021) | Examining the importance of bank risk exposures through interbank funding on bank efficiency levels | UK | Three-stage DEA with alliance between stage 1 and stage 2, and feedback variables | Net interbank fees, capital measured, Total interbank loans | Interbank borrowing, short-term wholesale funding and Long-term wholesale funding | 43 |
| (Omrani et al., 2022) | Analyzing the performance of branches of the Agribank in the West Azerbaijan province of Iran | Iran | Two-stage DEA-BWM model | Number of FTE employees, Total expenses, Number of ATM | Profit | 43 |
| Current Article | Evaluation of private bank branches | Iran | Dynamic DEA | Current facilities / exchange contracts, Total costs / operating expenses; | Total revenue / unrelated revenues | 531 |

- Clustering of the branches based on the results obtained from their efficiency.
- Performing a sensitivity analysis on the main variables to present a clear picture of each branch's status and upgrade each branch's status.
- Presenting a heuristic framework for dealing with big data and infeasible environments of DEAs.

Based on the provided review table of the literature, previous studies have used DEA for analyzing about 57 DMUs on average. In comparison, we used it for analyzing 531 DMUs which is 9.32 times greater than the existing problems in Table 1. The basic DEA models could not find feasible solutions in our case study. Therefore, we present a heuristic framework that can analyze a wide range of DMUs.

The main contributions of this paper are as follows:

1- Proposing a new heuristic framework for implementing DEA algorithms under a big data environment comprising many DMUs.
2- Employing dynamic DEA (DDEA) on a real case study including many DMUs equipped with a cluster-based hybrid method.
3- Presenting an efficient approach for rectifying infeasibility problems in DEA algorithms, so that the exposing probability of infeasible solutions reduces.
4- Creating the possibility of employing the proposed framework for every DEA algorithm.

The remainder of this paper is organized as follows. Section 2 provides a description of the studied problem, the mathematical model, and the relevant details of the developed framework. Section 3 implements the developed framework on a synthetic dataset and shows its capability and applicability. Section 4 summarizes the findings and suggests several paths for future research.

## 3. Methodology

The statistical population of this research comprises all 531 branches of a private bank in Iran, and according to the type of method used, there will be no need for sampling. The required data is provided to the researchers by the Studies and Innovation Center of the bank. However, to prepare the collected data for further analysis, it is necessary to perform the data cleansing process carefully.

Table 2: Introduction of branch financial ratios (final list of variables)

| **The ratio of financial statements and electronic transactions** | **Variable Names** |
|---|---|
| Current facilities/exchange contracts | $L_1$ |
| Total operating costs/expenses | $EX_1$ |
| Total non-operating expenses / expenses | $EX_2$ |
| Total revenue / common revenues | $IN_1$ |
| Total Income / separate revenues | $IN_2$ |

A comprehensive list of required indicators is collected based on the availability of data and the goals set by the bank managers. After an initial review of the variables from every aspect, including the review of the subject, consulting with financial experts, and justification of outcomes, the final variables used in the developed dynamic DEA are given in Table 2. All data were sorted by year. Each file contains separate sheets such as *Resources*, *Income*, *Facilities*, *Deferred*, and *Expense*, in which each data is sorted by the branch code. Please note that the first and second variables are *input*, the third variable is *undesirable link*, the fourth variable is *desirable link*, and the last variable is *output*.

Although DEA has been previously used in the banking industry, this method examines the system performance at a specific time. It does not provide an accurate understanding of the performance over time. In the present study, the dynamic SBM method is one of the most widely used models in this field (see Tone & Tsutsui, (2010)). This paper explains this model in detail. DEA offers many methods for measuring performance changes over time, such as a window analysis and the Malmquist Index. In this regard, it should be noted that these methods do not consider the relationships between different activities at two consecutive times and look at different times independently and separately. However, in a real business under long-term planning, considering the connections between different times helps to evaluate the system and gives a more accurate view. The Dynamic SBM technique considers the logical connections between the variables in consecutive time intervals.

These relationships are divided into four types: Desirable (good), undesirable (bad), free and fixed. The following is an overview of each of these types of communications:

- Desirable links such as net profit transferred to the profit carried forward, and net earned surplus carried to the next term.
- Undesirable links such as loss carried forward, bad debt, idle capital.
- Free and fixed links are optional and forced communications that may vary depending on the problem.

Suppose that we have the number *n* of DMUs whose information is specified in $T$ periods. Every DMU has a specific input and output at any given time. Also, both consecutive periods ($t$ and $t+1$) have specific correlations, shown in Fig. 1. Desirable links (optimal links) act as outputs. In this case, the value of the links should not be less than the observed value. Undesirable links act as input, where the value of the links should not exceed the observed value. The comparative lack of links in this section is called *inefficiency*.

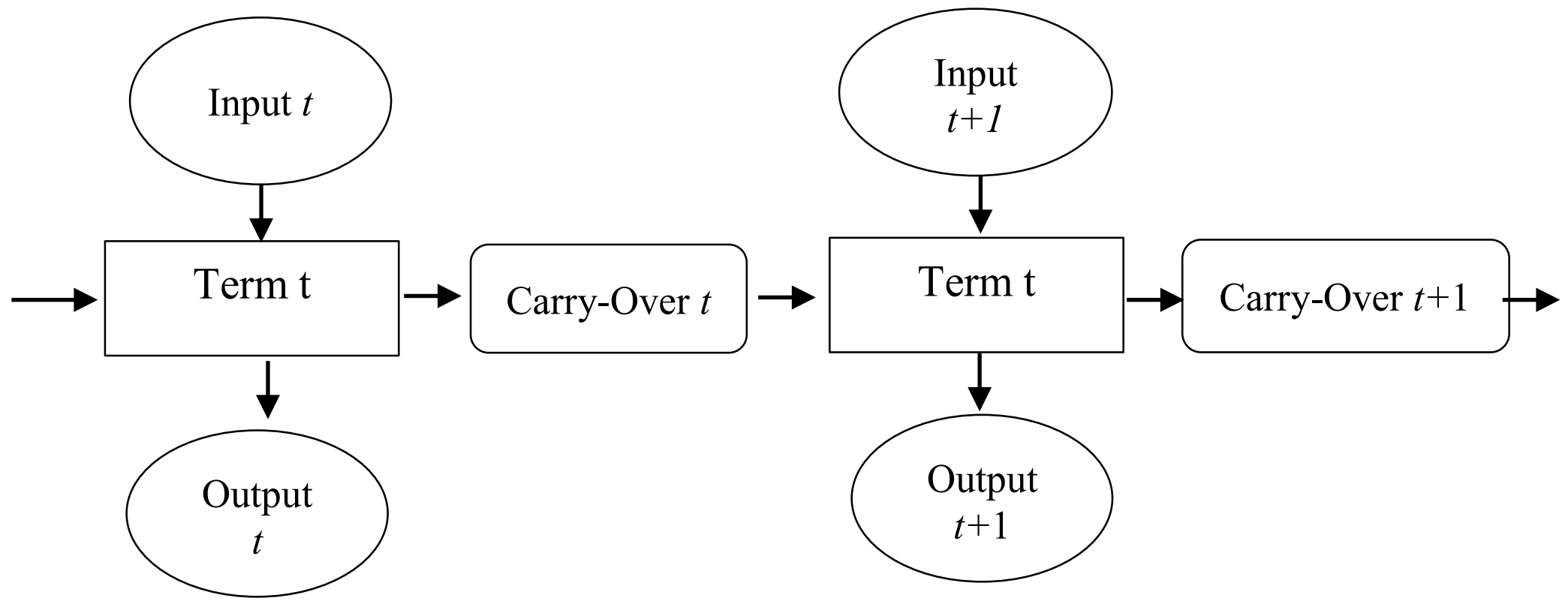

Fig. 1 Correlation between the time periods *t* and *t + 1* in a hypothetical DMU

Fig. 1 is drawn for a specific DMU. As can be seen, there are inputs and outputs at any given time,

and different times for a DMU are connected through the links. Sometimes a situation occurs, where the problem's nature requires only one of these items to be presented in the model. We first define the data in the dynamic DEA model, which are explained in Table 3. The '$z_{ijt}^{good}$' signifies good link values, where '$n_{bad}$' and '$n_{good}$' are the number of bad and good links, repsectively.

Table 3: Inputs and outputs of dynamic data envelopment analysis model

| Data | Index | Notes |
|---|---|---|
| $x_{ijt}$ | $i$ = 1,2,3,…,$m$; $j$=1,2,3,…,$n$; $t$=1,2,3,…,$T$ | Input |
| $y_{ijt}$ | $i$ = 1,2,3,…,$s$; $j$=1,2,3,…,$n$; $t$=1,2,3,…,$T$ | Output |
| $z_{ijt}^{good}$ | $i$ = 1,2,3,…,$n_{good}$; $j$=1,2,3,…,$n$ ;$t$=1,2,3,…,$T$ | Desirable link |
| $z_{ijt}^{bad}$ | $i$ = 1,2,3,…,$n_{bad}$; $j$=1,2,3,…,$n$; $t$=1,2,3,…,$T$ | Undesirable link |
| $\lambda_j^t$ | $j$=1,2,3,…,$n$; $t$=1,2,3,…,$T$ | - |

## 3.1 Mathematical Model

The dynamic DEA model under the SBM policy can now be defined as follows. In the proposed model, only desirable and undesirable links are considered.

$$\rho = Min\left[\frac{\frac{1}{T}\sum_{t=1}^{T} W^t\left[k - \frac{1}{m + nbad}\left(\sum_{i=1}^{m}\frac{k.s_{it}^-}{x_{i0t}} + \sum_{i=1}^{nbad}\frac{k.s_{it}^{bad}}{z_{i0t}^{bad}}\right)\right]}{\frac{1}{T}\sum_{t=1}^{T} W^t\left[k - \frac{1}{s + ngood}\left(\sum_{i=1}^{s}\frac{k.s_{it}^+}{y_{i0t}} + \sum_{i=1}^{ngood}\frac{k.s_{it}^{good}}{z_{i0t}^{good}}\right)\right]}\right] \tag{1}$$

$s.t.$

$$\sum_{j=1}^{n}\lambda_j^t x_{ijt} + s_{it}^- = x_{i0t};\ i\ =\ 1,2,3,\dots,m;\ t = 1,2,3,\dots,T \tag{2}$$

$$\sum_{j=1}^{n}\lambda_j^t y_{ijt} - s_{it}^+ = y_{i0t};\ i = 1,2,3,\dots,s;\ t = 1,2,3,\dots,T \tag{3}$$

$$\sum_{j=1}^{n}\lambda_j^t z_{ijt}^{good} - s_{it}^{good} = z_{i0t}^{good};\ i\ =\ 1,2,3,\dots,n_{\text{good}};\ t = 1,2,3,\dots,T \tag{4}$$

$$\sum_{j=1}^{n}\lambda_j^t z_{ijt}^{bad} + s_{it}^{bad} = z_{i0t}^{bad};\ i\ =\ 1,2,3,\dots,n_{\text{bad}};\ t = 1,2,3,\dots,T \tag{5}$$

$$\sum_{j=1}^{n}\lambda_j^t z_{ijt}^{bad} = \sum_{j=1}^{n}\lambda_j^{t+1} z_{ijt}^{bad}; t = 1,2,3\dots,T-1 \tag{6}$$

$$\sum_{j=1}^{n} \lambda_j^t z_{ijt}^{good} = \sum_{j=1}^{n} \lambda_j^{t+1} z_{ijt}^{good}; t = 1,2,3 \ldots, T-1 \qquad (7)$$

$$\sum_{j=1}^{n} \lambda_j^t = 1; \quad t = 1,2,3 \ldots, T \qquad (8)$$

$$\lambda_j^t, \lambda_j^{t+1}, s_{it}^{bad}, s_{it}^{good}, s_{it}^{+}, s_{it}^{-} \geq 0$$

where $s_{it}^{+}$ ،$s_{it}^{-}$ ،$s_{it}^{good}$ ،$s_{it}^{bad}$ و $s_{it}^{free}$ are slack variables which are described in Table 4.

Table 4: Description of the types of deficiency variables used in the main model

| **Variable Name** | **Index** | **Notes** |
|---|---|---|
| $\boldsymbol{s_{it}^{bad}}$ | $i = 1,2,3,\ldots,n_{bad}$; $t=1,2,3,\ldots,T$ | Undesirable surplus |
| $\boldsymbol{s_{it}^{good}}$ | $i = 1,2,3,\ldots,n_{good}$ ; $t=1,2,3,\ldots,T$ | Lack of optimal source |
| $\boldsymbol{s_{it}^{+}}$ | $i = 1,2,3,\ldots,s$; $t=1,2,3,\ldots,T$ | Lack of output |
| $\boldsymbol{s_{it}^{-}}$ | $i = 1,2,3,\ldots,m$; $t=1,2,3,\ldots,T$ | Input surplus |
| $\boldsymbol{\lambda_j^t}$ | $j = 1,2,3,\ldots,n$; $t=1,2,3,\ldots,T$ | Link coefficient |

The primary goal of this method (Eq. 1) is to separate efficient and inefficient units and rank them. Although the proposed model is nonlinear, it can be modified to become a linear model. This is achieved by setting the denominator of the fraction in the objective function equal to one and adding it to the constraints. The objective function states that the maximum output limit is obtained for the smallest input values. As can be seen, the two variables of desirable and undesirable links have also been added to the Dynamic SBM model, Eqs. (4) and (5). Undesirable links are input, where we tend to reduce them. On the other hand, the desirable links are of the output type, and we tend to increase their values. The inputs should not exceed a specific limit. Therefore, we set restrictions such as $x_{i0t}$ (Eq. 2) and $z_{i0t}^{bad}$ (Eq. 5) on the inputs and undesirable variables, respectively. This is opposite to the output constraints and the constraints on the desirable variables. Since we want these variables not to be less than a specific limit, their constraints are defined to be greater than equal (Equations 3 and 4). We connect consecutive years by the dynamic nature of the model (in which the time factor is involved). Equations (6) and (7) are used for this purpose. The desirable and undesirable links are linked for two consecutive years. With such an objective function, these constraints must be resolved separately for each unit. Eq. (8) related to the variable return-to-scale assumption.

### 3.2 DEA framework

As shown in Table 1, most papers have utilized DEA for comparing a low number of DMUs, resulting in a low number of constraints in the DEA model, increasing the feasible solution area. Indeed, as the number of DMUs increases, the feasible solution space becomes smaller, and the probability of reaching a feasible answer increases. However, this paper fed a wide range of DMUs (531 DMUs) to the DEA model, which could not reach feasible solutions. Accordingly, a novel heuristic framework deals with infeasible solutions for implementing DEA algorithms for many DMUs.

According to the proposed heuristic framework, each DMU is considered first, and the rest of the DMUs are randomly categorized into $p$ equal classes. Each time, this opted DMU is added to each cluster, and then the efficiency of this selected DMU is calculated. This process should be conducted for all of the clusters. By doing so, $p$ numbers, as efficiency values, are obtained. Then, it suffices to take the average of these $p$ numbers, signifying the estimated efficiency of this opted DMU. This process should be carried out for the 531 DUMs considered in this research. The proposed framework can be seen in Fig. 2. The goal of taking the average of the obtained numbers is to increase the robustness of the outputs.

### 3.3. Clustering Process

Clustering is unsupervised learning in which samples are divided into groups whose members are similar to each other. We utilize the k-means algorithm to cluster the branches. Among the essential features of a K-Means clustering, we can refer to

1. the algorithm's performance depends on the primary centers,
2. trapping into a local optimum,
3. the convexity of the shape of clusters, and
4. a practical algorithm for large volumes of data.

The K-means algorithm considers an $m \times n$ matrix as $m$ data with $n$ properties; If $k$ is 2, it divides the data into two categories, where the members in the same category have properties that are closer to each other. Consider the set of vectors $a_1, a_2, \dots, a_n$. The K-means algorithm is used to find the clusters $\pi_1, \pi_2, \dots, \pi_k$ using the following function (Eq. 9):

$$D\left(\{\pi_c\}_{c=1}^{k}\right) = \sum_{c=1}^{k} \sum_{a_{i\in\pi_i}} \|a_i - m_c\|^2 \quad (9)$$

where $m_c = \sum_{a_i \in \pi_i} a_i / |\pi_c|$ indicates the average of the $c$'th cluster, where cluster $c$ is represented by $\pi_c$.

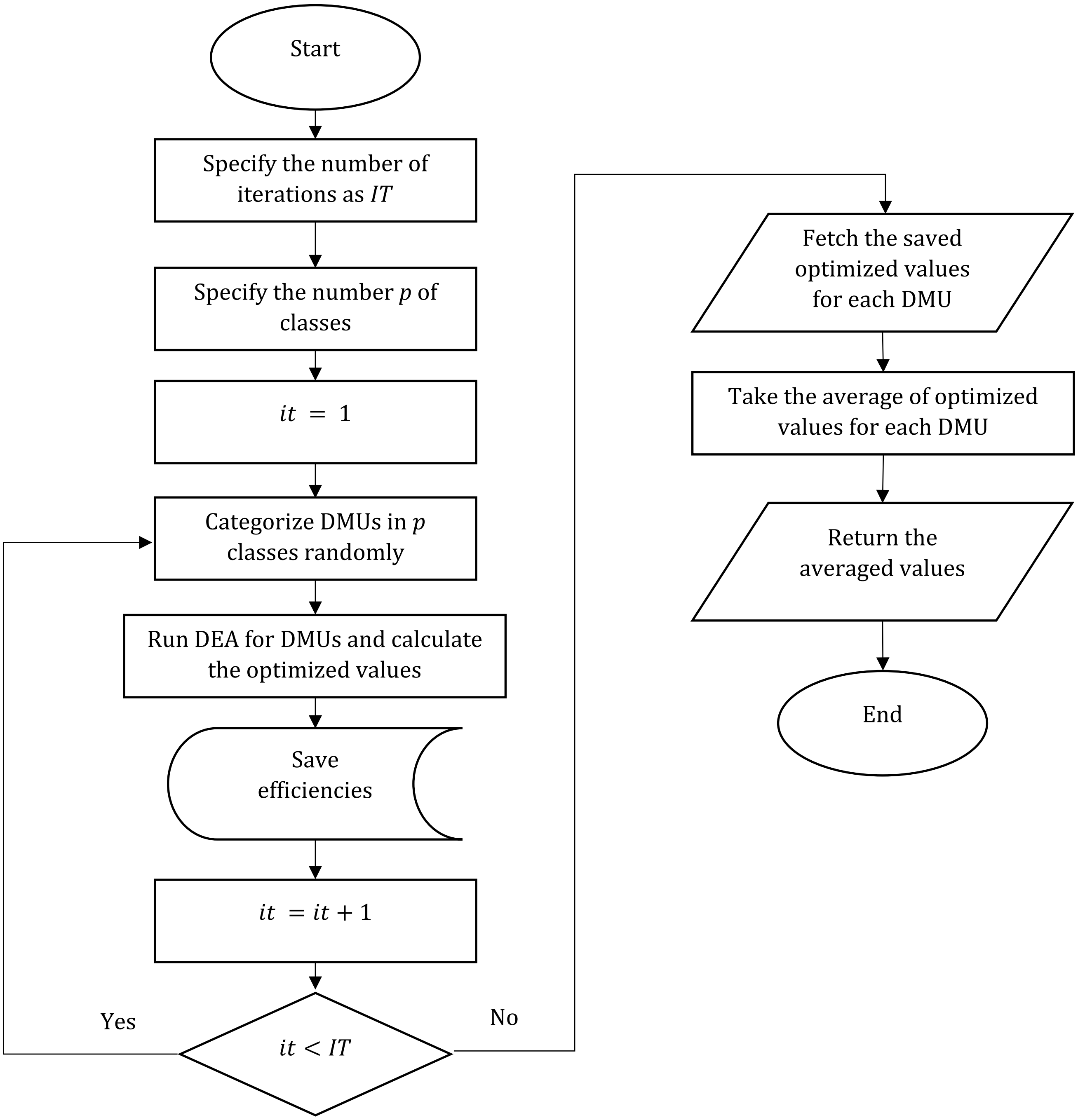

Fig. 2 Proposed heuristic framework of DEA

# 4. Experimental Results

In this section, the results of the proposed method are presented. We use Python 3.7 on an Intel Core i7 processor (2.7 GHz) with 12GB memory to run the experiments. After running the DDEA algorithm as described in the previous section, the top 20 branches of the private bank in Iran between 2017 and 2019 are listed in Table 5 with the current grade of the branch, the efficiency value, and the branch code. As can be seen, most of these branches are first-class, and there are only eight branches of other degrees, of which six are second-class.

Table 5: Details of the top 20 branches of the case study during the period from 2017 to 2019 using the DDEA method

| Branch Type | Efficiency | Branch Code | Rank |
|---|---|---|---|
| Special | 2.7845 | 3610 | 1 |
| Premier | 2.5389 | 1867 | 2 |
| Premier | 2.5379 | 1887 | 3 |
| Premier | 2.5166 | 1857 | 4 |
| 1 | 2.5080 | 3410 | 5 |
| 2 | 2.4219 | 3511 | 6 |
| 2 | 2.3462 | 2570 | 7 |
| 2 | 2.2124 | 4614 | 8 |
| 4 | 2.1660 | 1202 | 9 |
| 2 | 2.0125 | 1225 | 10 |
| 2 | 2.0065 | 2061 | 11 |
| Premier | 1.9826 | 1839 | 12 |
| Premier | 1.9746 | 1818 | 13 |
| 1 | 1.8744 | 1110 | 14 |
| Premier | 1.8469 | 1870 | 15 |
| 3 | 1.7964 | 1631 | 16 |
| Special | 1.7895 | 4410 | 17 |
| Premier | 1.7260 | 1885 | 18 |
| 2 | 1.6880 | 3310 | 19 |
| 1 | 1.6790 | 2510 | 20 |

Table 6: Descriptive statistics related to the variables used for the three years

| The Ratio of Financial Statements and Transactions | Variable Name | Min | Max | Variance | Year |
|---|---|---|---|---|---|
| Current facilities/exchange contracts | $L_1$ | 0 | 0.0520 | 0.0042 | 2107 |
| | | 0 | 0.519 | 0.00416 | 2018 |
| | | 0 | 0.0517 | 0.0042 | 2019 |
| Total operating costs/expenses | $EX_1$ | 0 | 0.0532 | 0.0037 | 2107 |
| | | 0 | 0.522 | 0.0035 | 2018 |
| | | 0 | 0.0512 | 0.0037 | 2019 |
| Total non-operating expenses / expenses | $EX_2$ | 0 | 0.079 | 0.0039 | 2017 |
| | | 0 | 0.0820 | 0.004 | 2018 |
| | | 0 | 0.0810 | 0.0039 | 2019 |
| Total revenue / common revenues | $IN_1$ | 0 | 0.0799 | 0.0048 | 2017 |
| | | 0 | 0.087 | 0.0041 | 2018 |
| | | 0 | 0.087 | 0.0039 | 2019 |
| Total Income / Unrequited Income | $IN_2$ | 0 | 0.0724 | 0.0038 | 2017 |
| | | 0 | 0.0713 | 0.0044 | 2018 |
| | | 0 | 0.0753 | 0.0040 | 2019 |

The information about the variables is shown in Table 6 in the order of the years. The range of variables is almost the same in consecutive years, which indicates that there are not many outliers in the database that may deviate the final output in a certain direction.

### 4.1 Clustering Implementation

This study uses the silhouette criterion to evaluate different clustering solutions. Fig. 3 shows different silhouette values for a number of clusters between 2 and 12. As can be seen, the silhouette for 8 clusters has reached its maximum value of 78.21%, and it is equal to 76.48% for 7 clusters.

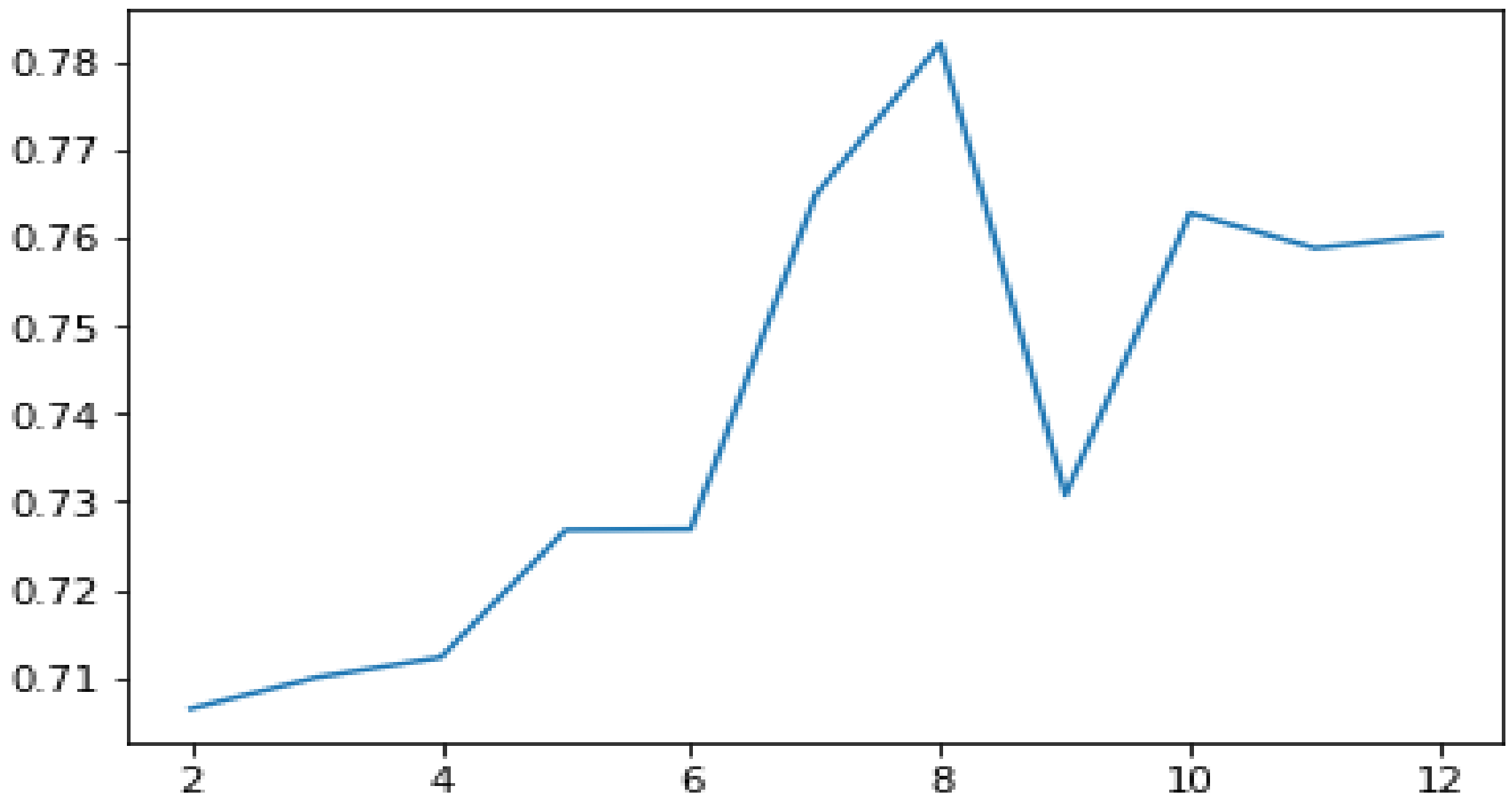

Fig. 3 The standard value of silhouette for different clusters

Based on the Silhouette criterion, seven clusters were disclosed. There are 7 degrees, including *privileged*, *special*, *1*, *2*, *3*, *4*, and *5*, in the classification of the bank branches.

The best silhouette values were calculated for clustering algorithms such as hierarchical clustering Ward, K-means, DBSCAN, and OPTICS. The best silhouette values of the K-means method are larger than the amount of the silhouette obtained from the other methods. Therefore, in this research, the K-means method has been used. After ranking the bank branches in the three years, the clustering operation was performed using the K-means method with the number of 7 clusters on the output obtained from the ranking. In this clustering, groups 1 and 2 can be considered as equivalent to special and privileged branches, and the rest of the clusters represent the branches of grades 1, 2, 3, 4, and 5, respectively. Fig. 4 shows the percentage of dispersion of each of these types.

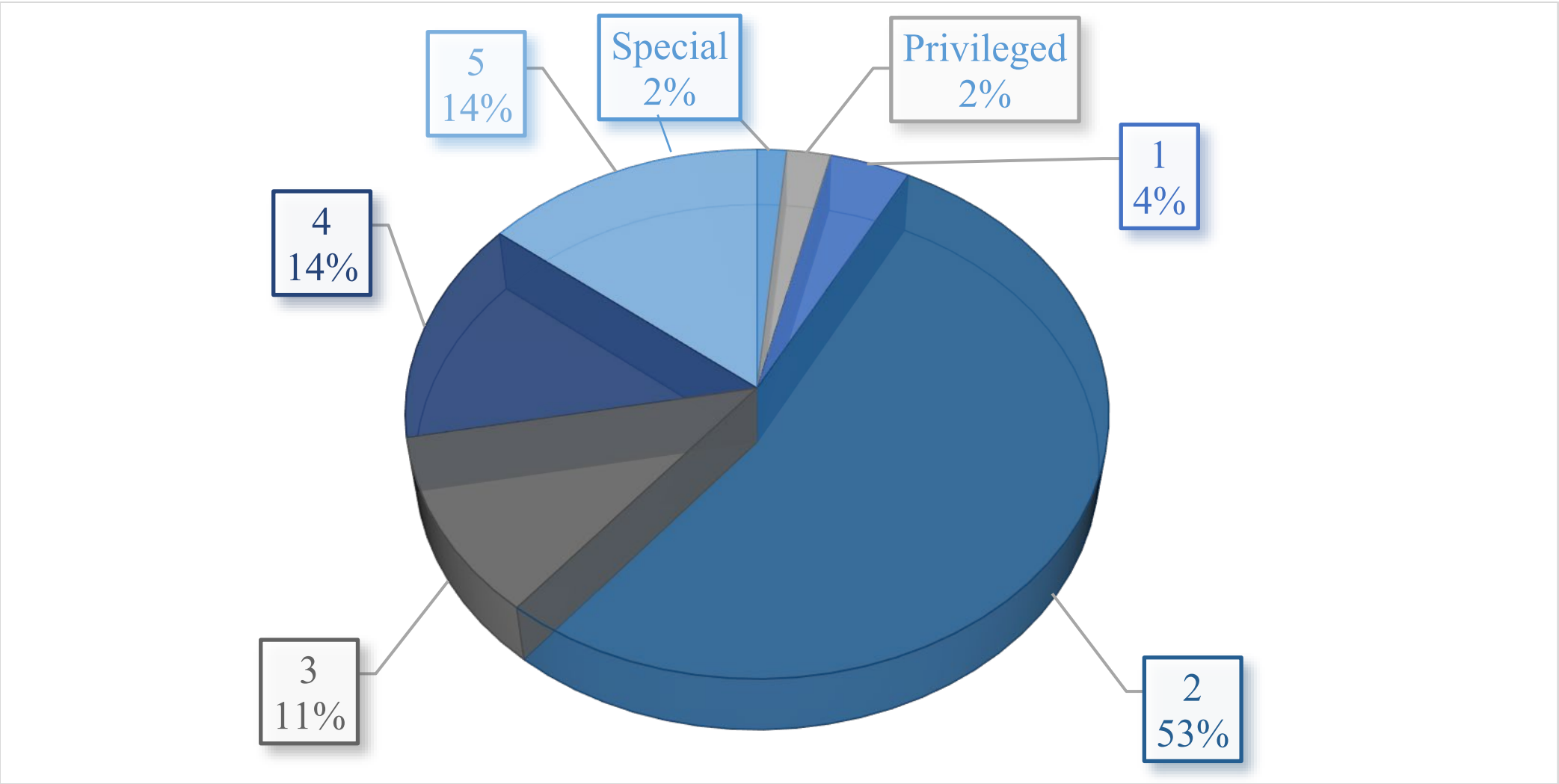

Fig. 4 The share of different types of bank branches in the clustering

The clustering results of the top 20 branches of the bank during the three years are given in Table 7.

Table 7: Clustering results of the top 20 branches of the bank during the years from 2017 to 2019

| Cluster # | Performance value | Branch Code | Rank |
|---|---|---|---|
| **1** | 2.7845 | 3610 | 1 |
| **1** | 2.5389 | 1867 | 2 |
| **1** | 2.5379 | 1887 | 3 |
| **1** | 2.5166 | 1857 | 4 |
| **1** | 2.5080 | 3410 | 5 |
| **1** | 2.4219 | 3511 | 6 |
| **1** | 2.3462 | 2570 | 7 |
| **1** | 2.2124 | 4614 | 8 |
| **2** | 2.1660 | 1202 | 9 |
| **2** | 2.0125 | 1225 | 10 |
| **2** | 2.0065 | 2061 | 11 |
| **2** | 1.9826 | 1839 | 12 |
| **2** | 1.9746 | 1818 | 13 |
| **2** | 1.8744 | 1110 | 14 |
| **2** | 1.8469 | 1870 | 15 |
| **2** | 1.7964 | 1631 | 16 |
| **2** | 1.7895 | 4410 | 17 |
| **2** | 1.7260 | 1885 | 18 |
| **2** | 1.6880 | 3310 | 19 |
| **2** | 1.6790 | 2510 | 20 |

In the clustering, the center values of the clusters are used to compare the clusters. The use of these centers gives a proper representation of the distance between two clusters or the similarity of the two clusters. The efficiency of the cluster centers is shown in Table 8.

Table 8: The efficiency of the cluster centers

| **Cluster Number** | **Cluster Type** | **Efficiency** |
|---|---|---|
| **1** | Special | 2.48 |
| **2** | Privileged | 1.878 |
| **3** | 1st Rank | 1.311 |
| **4** | 2nd Rank | 1.006 |
| **5** | 3rd Rank | 0.781 |
| **6** | 4th Rank | 0.473 |
| **7** | 5th Rank | 0.033 |

According to Table 8, the efficiency of "special" branches is, on average, equal to 2.48, considered as "super-efficient" branches. It is important to note that in the SBM method, the efficiency of a branch can reach more than one, which indicates that the branch is "super-efficient".

## 4.2 Sensitivity Analysis

A sensitivity analysis shows the amount of change in each of the five key variables studied, intending to upgrade the level of each branch to a higher level. In other words, conducting such a sensitivity analysis allows the bank managers to assess the current state of the branch over three years and take into account the average performance of each branch in the medium term, which is naturally a more accurate judgement of its status. It also helps to mitigate the negative effect of weaknesses or temporary growth of branches and to increase awareness of branch managers on the amount of the necessary change in each index by each branch. The awareness of the tolerance of these changes gives the bank managers a better view of the allowable change in each index.

Table 9 shows the results of the sensitivity analysis of 12 branches located in the second cluster (privileged branches). Based on these results, the change needed in each branch to upgrade that branch to a higher cluster is given. Due to the simultaneous change of the considered variables, the values listed in Table 9 are calculated based on the difference in the number of each branch's variables compared to the worst branch in the higher cluster.

Table 9: Results of the sensitivity analysis for the second cluster branches (Privileged branches)

| **Non-operating Costs** | **Common Income** | **Separate Income** | **Operating Costs** | **Facilities** | **Branch Code** |
|---|---|---|---|---|---|
| **-3332.33** | No Change | 6265.667 | No Change | No Change | 1110 |
| **No Change** | No Change | 48904.67 | No Change | No Change | 1202 |
| **-1961.67** | 20978.66 | 30688.67 | No Change | No Change | 1225 |
| **No Change** | 22063.33 | 5294 | No Change | 22974 | 1631 |
| **-19686** | No Change | No Change | -62295.7 | No Change | 1818 |
| **-60950.3** | No Change | No Change | -974662 | No Change | 1839 |
| **-5788.67** | No Change | No Change | -75570.3 | No Change | 1870 |
| **-50661.7** | No Change | No Change | -26264.3 | No Change | 1885 |
| **No Change** | 30440.66 | 69275.67 | No Change | 47038.33 | 2061 |
| **-12622** | No Change | No Change | -25587.7 | No Change | 2510 |
| **-2813.33** | 27814.66 | No Change | No Change | No Change | 3310 |
| **-155713** | No Change | No Change | -798827 | No Change | 4410 |

For further clarification, we describe the details of one of the branches, branch 1110, located in the second cluster. We want to know what changes can upgrade this branch to a higher cluster, i.e., cluster 1 (special branches). For this purpose, we look up the information related to the same row in the columns related to the variables of cluster 2. According to Table 9, this branch should increase its separate income by about 6266 units and reduce its non-operating expenses by about 3332 units simultaneously. While improving the two mentioned variables, the other variables should not be diminished. No discount should occur so that this branch can be upgraded to the first cluster (if the indicators related to the worst branch of cluster 1 remain constant).

Table 10 shows the results of the sensitivity analysis related to the branches of the third cluster (Grade 1 branches). Based on this sensitivity analysis, it is possible to understand how much change is needed to turn a grade 1 branch into a premium branch. As can be seen, there is no need to change some of these indicators.

Table 10: Results of the sensitivity analysis for the third cluster branches (Grade 1 branches)

| Non-operating Costs | Common Income | Separate Income | Operating Costs | Facilities | Branch Code |
|---|---|---|---|---|---|
| **No Change** | 57071.67 | 107587 | No Change | 1687068 | 1128 |
| **No Change** | 21787.67 | No Change | No Change | 1523104 | 1313 |
| **No Change** | 53308 | 77444.3 | No Change | 1687718 | 1440 |
| **No Change** | 49338.67 | 59115.3 | No Change | 1670713 | 1448 |
| **-7416.33** | 46334.67 | No Change | -11796.3 | 1209070 | 1861 |
| **No Change** | 36147.67 | No Change | -39408.3 | 1417276 | 1862 |
| **No Change** | No Change | 78671.3 | No Change | 1569883 | 1866 |
| **-53387** | No Change | No Change | -270581 | No Change | 1886 |
| **No Change** | 56365.67 | 118305.3 | No Change | 1698959 | 2047 |
| **No Change** | 36198.67 | 85478.3 | No Change | 1669256 | 2115 |
| **No Change** | 7905 | 68022.97 | No Change | 1522511 | 2310 |
| **No Change** | 64203.67 | 173368.6 | No Change | 1738157 | 2611 |
| **No Change** | 48579.67 | No Change | -15437.7 | 1532845 | 2734 |
| **No Change** | 53547.33 | 62155.3 | No Change | 1422109 | 3110 |
| **No Change** | 55061 | 159659.3 | No Change | 1734663 | 3318 |
| **No Change** | 56431 | No Change | No Change | 1619914 | 3411 |
| **No Change** | 31804.67 | 15928.97 | No Change | 1625130 | 3514 |
| **-8047.33** | No Change | No Change | -25805.7 | No Change | 3813 |
| **No Change** | 45289.67 | 63090.63 | No Change | 1109168 | 3822 |
| **No Change** | 53035.67 | 188829.3 | No Change | 1766943 | 4020 |
| **No Change** | 29063.33 | 83175.63 | No Change | 1625689 | 4617 |
| **No Change** | No Change | No Change | -10012 | 1188718 | 4622 |

Similarly, Tables A.1-A.4 illustrate the sensitivity analysis of Grade 2 (4th cluster), Grade 3 (5th Cluster), Grade 4 (6th Cluster), and Grade 5 (7th Cluster) branches, respectively. Based on this sensitivity analysis, it is possible to determine how many changes are needed to turn a lower-grade branch into a higher-grade branch. Similar to the previous tables, some sections do not need to be changed. This means that these indexes are good enough, considering the values of the change of the other variables as needed. As mentioned, this main philosophy of the sensitivity analysis in this study is to compare each cluster's worst branch. Thus, if a branch wants to be upgraded from a 4th-degree cluster (worst cluster) to a 3rd-degree cluster (better degree), it must do better than or equal to at least one of the worst members of the 3rd cluster. Let us explain this by an example. Suppose that the worst branch of cluster 3 has a cost of 10 units and the branch we consider, type 4, has a cost of 15. Suppose that this variable is the only criterion for comparing the two branches. In such a situation, if the branch under consideration can at least perform the same as the worst branch of cluster 3 (i.e., reduce its costs to 10 units), it must reduce its costs by at least 10-5 = 5 units. Finally, it should be noted that the present study used dynamic data envelopment analysis to evaluate the 531 branches. However, no article in the existing literature considers more than 50 branches. One of the reasons for the lack of research could be the high computational volume and complexity of the algorithm used, which prevents a large number of branches from being considered.

### 4.3 Comparison of the results of the DEA and Dynamic DEA methods

It is not common to compare the results of the classic SBM and DEA models, as the range of variation between the two methods is different. In the DEA, the performance values are generally between zero and one, while in SBM-based models, the performance values include numbers more significant than one. Nevertheless, we can perform a relative comparison of these two models. It should also be noted that comparing DEA with dynamic DEA (DDEA) would be reasonable and common in the DDEA models with performance values between zero and one. The values obtained in the DEA method for the branches of the studied bank are consistent with the results obtained from the DDEA method. A part of this comparison is illustrated in Table 11.

Table 11: Comparison of DDEA and DEA

| Branch Code | Efficiency 2017 | Efficiency 2018 | Efficiency 2019 | DDEA Efficiency |
|---|---|---|---|---|
| **3610** | 1 | 1 | 1 | 2.7854 |
| **1867** | 1 | 1 | 1 | 2.5389 |
| **1887** | 1 | 1 | 1 | 2.5378 |
| **1857** | 1 | 1 | 0.9916 | 2.5160 |
| **3410** | 1 | 1 | 0.9880 | 2.5079 |

As can be seen, the first five branches in all three time periods have numbers very close to one. This is a testament to the accuracy of the performance obtained from the dynamic data envelopment analysis method.

## 5. Conclusions and Future Works

Evaluating the efficiency of organizations and companies or evaluating different branches of an organization is one of the challenging issues that the management has faced. To date, various research has been conducted on this issue. Among the many methods used in the evaluation of bank branches, researchers have become more inclined to non-parametric methods in recent years and constitute a major part of the evaluation of organizational branches. Among these, one of the most widely used non-parametric methods are data envelopment analysis techniques. Data envelopment analysis can be divided into static and dynamic parts, where dynamic DEA leads to results closer to the real world because the real environment changes over time and is not static. Excluding time can increase the possibility of errors in the output of the work and as a result, a wrong analysis of the output increases. After ranking the branches, a clustering of the branches was performed using the K-means method.

The silhouette index was used to find the optimal number of clusters, which was the best value for the 8 clusters. However, due to the very small difference in the index value for 7 and 8 clusters and the classification of the bank branches to 7 degrees and to equalize and make the analysis more uniform, we use 7 clusters for clustering the operations. In this clustering, groups 1 and 2 are equivalent with special and privileged branches, while the rest of the clusters show branches of grades 1, 2, 3, 4, and 5, respectively. A comprehensive sensitivity analysis was performed on the results. The sensitivity results show the rank of each branch and the rate of change in each of the five key indicators studied to raise the level of each branch to a higher level. The sensitivity analysis allows the bank managers to recognize the amount of necessary change in each index by each branch. The awareness of the rate of change of each indicator to lower the level of each branch to a lower degree can also be considered in the plans of the bank's general managers. This helps the managers to prevent a violation of the set values and prevent the degradation of the branch.

In the last part of the experiments, the classical DEA method was compared with the DDEA method. The results of this comparison showed how the output obtained from the DDEA method is in full compliance with the three-year analysis of the bank branches through the classical DEA method. The present study used dynamic data envelopment analysis to evaluate 531 branches. However, it can be said that almost no article in the existing literature considers more than 50 branches. Regarding the limitations of the problem under consideration, the "complexity" of the DDEA method can be considered as the most important challenge ahead. In other words, considering a larger number of variables (for example, ten variables) will most likely lead to an infeasible answer. The following future works can be a step towards a complete research:

- Providing a smart and user-friendly dashboard for bank experts to take advantage of the method implemented in this research.
- Comparing with other modern methods; this is one of the cases that can always be considered to complete the impact of any research.
- Combining/hybridizing different methods; this should be done considering the strengths and weaknesses of different methods.
- Utilizing big data-based methods; since the data volume is large for performing DEA dynamic analysis, and it takes a lot of time to execute the relevant code, it can be facilitated to implement this method. In this regard, using parallel algorithms and discrete processing seems to be necessary.
- Examining the relationships between variables; Studies show that most articles use domain experts to determine their study variables. Therefore, providing an unsupervised method to identify input and output variables can save time. Also, measuring the correlation coefficient between these variables and such methods can improve the quality of solving such a problem.